\documentclass[11pt]{article}
\usepackage{amsmath,amssymb}

\def\N{\bf \mbox{I\hspace{-.15em}N}}

\def\R{{\bf \mbox{I\hspace{-.20em}R}}}

\def\bkR{{\bf \mbox{I\hspace{-.20em}R}}}
\def\C{C^{\infty}(M, {\bf \mbox{I\hspace{-.20em}R}})}

\def\CR{C^{\infty}(M\times {\bf \mbox{I\hspace{-.20em}R}}, {\bf \mbox{I\hspace{-.20em}R}})}

\newtheorem{definition}{Definition}[section]

\newtheorem{proposition}[definition]{Proposition}
\newtheorem{theorem}[definition]{Theorem}
\newtheorem{remark}[definition]{Remark}

\newtheorem{examples}[definition]{Examples}

\newenvironment{proof}{\noindent{\textbf{Proof}.}}{\hfill $\square$}

\begin{document}

\title{On quasi-Jacobi and Jacobi-quasi bialgebroids}
\author{J. M. Nunes da Costa \\Departamento de Matem\'atica \\
Universidade de Coimbra\\Apartado 3008
\\3001-454 Coimbra - Portugal\\ {\small  E-mail: jmcosta@mat.uc.pt}
\and F. Petalidou\\Faculty of Sciences and Technology
\\University of Peloponnese \\22100 Tripoli - Greece\\{\small E-mail: petalido@uop.gr}}

\date{}
\maketitle

\begin{abstract}
We study quasi-Jacobi and Jacobi-quasi bialgebroids and their
relationships with twisted Jacobi and quasi Jacobi manifolds. We
show that we can construct quasi-Lie bialgebroids from
quasi-Jacobi bialgebroids, and conversely, and also that the
structures induced on their base manifolds are related via a
``quasi Poissonization".
\end{abstract}

\vspace{3mm} \noindent {\bf{Keywords: }}{Quasi-Jacobi bialgebroid,
Jacobi-quasi bialgebroid, quasi-Jacobi bialgebra, twisted Jacobi
manifold, quasi Jacobi manifold.}

\vspace{3mm} \noindent {\bf{A.M.S. classification (2000):}} 17B62;
17B66; 53D10; 53D17.

\section{Introduction}
The notion of {\em quasi-Lie bialgebroid} was introduced in
\cite{ro}. It is a structure on a pair $(A,A^*)$ of vector
bundles, in duality, over a differentiable manifold $M$ that is
defined by a Lie algebroid structure on $A^*$, a skew-symmetric
bracket on the space of smooth sections of $A$ and a bundle map $a
: A \to TM$, satisfying some compatibility conditions. These
conditions are expressed in terms of a section $\varphi$ of
$\bigwedge^3A^*$, which turns out to be an obstruction to the Lie
bialgebroid structure on $(A,A^*)$. A quasi-Lie bialgebroid will
be denoted by $(A,A^*,\varphi)$. In the case where $A$ is a Lie
algebroid, its dual vector bundle $A^*$ is equipped with a
skew-symmetric bracket on its space of smooth sections and a
bundle map $a_* :A^* \to TM$ and the compatibility conditions are
expressed in terms of a section $Q$ of $\bigwedge^3A$, the triple
$(A,A^*, Q)$ is called a \emph{Lie-quasi bialgebroid} \cite{ks3}.
When $\varphi=0$ and $Q=0$, quasi-Lie and Lie-quasi bialgebroids
are just Lie bialgebroids. We note that, while the dual of a Lie
bialgebroid is itself a Lie bialgebroid, the dual of a quasi-Lie
bialgebroid is a Lie-quasi bialgebroid, and conversely \cite{ks3}.
The quasi-Lie and Lie-quasi bialgebroids are particular cases of
proto-bialgebroids \cite{ks3}. As in the case of a Lie
bialgebroid, the doubles $A\oplus A^*$ of a quasi-Lie and of a
Lie-quasi bialgebroid are endowed with a Courant algebroid
structure \cite{ro}, \cite{ks3}.

It was shown in \cite{ro1} that the theory of quasi-Lie
bialgebroids is the natural framework in which we can treat
\emph{twisted Poisson manifolds}. These structures were introduced
in \cite{sw}, under the name of Poisson manifolds with a closed
$3$-form background, motivated by problems of string theory
\cite{p} and of topological field theory \cite{kl}.

The notion of \emph{Jacobi bialgebroid} and the equivalent one of
\emph{generalized Lie bialgebroid} were introduced, respectively,
in \cite{gm1} and \cite{im1}, in such a way that a Jacobi
bialgebroid is canonically associated to a Jacobi manifold and
conversely. A Jacobi bialgebroid over $M$ is a pair
$((A,\phi),(A^*,W))$ of Lie algebroids over $M$, in duality,
endowed with $1$-cocycles $\phi\in \Gamma(A^*)$ and $W\in
\Gamma(A)$ in their Lie algebroid cohomology complexes with
trivial coefficients, respectively, that satisfy a compatibility
condition. Also, its double $(A\oplus A^*,\phi +W)$ is endowed
with a Courant-Jacobi algebroid structure \cite{gm2}, \cite{jj}.

In order to adapt to the framework of Jacobi manifolds the
concepts of twisted Poisson manifold and quasi-Lie bialgebroid, we
have recently introduced in \cite{jf-tj} the notions of
\emph{twisted Jacobi manifold} and \emph{quasi-Jacobi
bialgebroid}.  The purpose of the present paper is to develop the
theory of quasi-Jacobi bialgebroids, as well as of its dual
concept of \emph{Jacobi-quasi bialgebroids}, and to establish a
very close relationship between quasi-Jacobi and quasi-Lie
bialgebroids.

\vspace{.2cm}

The paper contains four sections, besides the Introduction, and
one Appendix (section 5). In section 2 we recall the definition of
quasi-Jacobi bialgebroid, we present some basic results
established in \cite{jf-tj}, we develop the examples of
quasi-Jacobi and Jacobi-quasi bialgebroids associated to twisted
Jacobi manifolds and to quasi Jacobi manifolds, and, finally, we
study the triangular quasi-Jacobi bialgebroids. Section 3 is
devoted to the study of the structures induced on the base
manifolds of quasi-Jacobi and Jacobi-quasi bialgebroids. Several
examples are presented. In section 4 we establish a one to one
correspondence between quasi-Jacobi bialgebroid structures
$((A,\phi),(A^{\ast},W),\varphi)$ over a manifold $M$ and
quasi-Lie bialgebroid structures $(\tilde{A},\tilde{A}^{\ast},
\tilde{\varphi})$ over $\tilde{M}=M\times \R$. Also, we prove that
the  structure induced on $\tilde{M}=M\times \R$ by
$(\tilde{A},\tilde{A}^{\ast}, \tilde{\varphi})$ is the ``quasi
Poissonization" of the structure induced on $M$ by
$((A,\phi),(A^{\ast},W),\varphi)$. The dual version of these
results is also presented. Finally, in the Appendix, we define the
action of a Lie algebroid with $1$-cocycle on a differentiable
manifold, a concept that is used in the paper.

\bigskip

\noindent {\bf Notation:} If $(A,\phi)$ is a Lie algebroid with
$1$-cocycle $\phi$, we denote by $d^\phi$ the differential
operator $d$ of $A$ modified by $\phi$, i.e., $d^\phi\alpha= d
\alpha + \phi \wedge \alpha$, for any $\alpha \in
\Gamma(\bigwedge^kA^*)$. Moreover, we denote by $\delta$ the usual
de Rham differential operator on a manifold $M$ and by ${\rm d}$
the differential operator of the Lie algebroid $TM \times \R$,
${\rm d}(\alpha, \beta)=(\delta \alpha, -\delta \beta)$, for
$(\alpha, \beta) \in \Gamma(\bigwedge^k (T^*M \times \R))\equiv
\Gamma(\bigwedge^k T^*M) \times \Gamma(\bigwedge^{k-1} T^*M)$. We
also consider the identification $\Gamma(\bigwedge^k (TM \times
\R))\equiv \Gamma(\bigwedge^k TM) \times \Gamma(\bigwedge^{k-1}
TM)$.

\section{Quasi-Jacobi and Jacobi-quasi bialgebroids}
Let $((A,\phi),(A^*,W))$ be a pair of dual vector bundles over a
differentiable manifold $M$ endowed with a $1$-form $\phi$ and
$W$, respectively, and $\varphi$ a $3$-form of $A$.

\begin{definition} \label{def1.1}
A {\em quasi-Jacobi bialgebroid} structure on
$((A,\phi),(A^*,W),\varphi)$ consists of a Lie algebroid structure
with $1$-cocycle $([\cdot,\cdot]_*, a_*,W)$ on $A^*$, a bundle map
$a : A \to TM$ and a skew-symmetric operation $[\cdot, \cdot]$ on
$\Gamma(A)$ satisfying, for all $X,Y,Z \in \Gamma(A)$ and $f \in
\C$, the following conditions:
\begin{enumerate}
\item [1)] $[X,fY]=f[X,Y]+(a(X)f) Y$; \item [2)]
$a([X,Y])=[a(X),a(Y)]-a_* (\varphi(X,Y,\cdot))$; \item [3)]
$[[X,Y],Z]+c.p. = - d_*^W(\varphi(X,Y,Z)) -((i_{\varphi(X,Y,
\cdot)}d_*^W Z)+ c.p.)$; \item [4)] $d\phi -
\varphi(W,\cdot,\cdot)=0$, where $d$ is the quasi-differential
operator on $\Gamma(\bigwedge A^*)$ determined by the structure
$([\cdot, \cdot],a)$ on $A$; \item [5)] $d^\phi\varphi=0$,
\noindent $d^\phi$ being the quasi-differential operator modified
by $\phi$; \item [6)] $d_{*}^{W} [P,Q]^{\phi}= [d_{*}^{W}
P,Q]^{\phi} + (-1)^{p+1}[P, d_{*}^{W}Q]^{\phi}$, where
$[\cdot,\cdot]^{\phi}$ is the bracket on $\Gamma(\bigwedge A)$
modified by $\phi$, $P \in \Gamma (\bigwedge^{p}A)$ and $Q \in
\Gamma (\bigwedge A)$.
\end{enumerate}
\end{definition}

As in the case of quasi-Lie and Lie-quasi bialgebroids, by
interchanging the roles of $(A,\phi)$ and $(A^*,W)$ in the above
definition, we obtain the notion of \emph{Jacobi-quasi
bialgebroid} over a differentiable manifold $M$. We have: \emph{If
$((A,\phi),(A^*,W),\varphi)$ is a quasi-Jacobi bialgebroid over a
differentiable manifold $M$, then $((A^*,W),(A,\phi),\varphi)$ is
a Jacobi-quasi bialgebroid over $M$, and conversely.}

In the case where both $1$-cocycles $\phi$ and $W$ are zero, we
recover, from Definition  \ref{def1.1}, the notion of quasi-Lie
bialgebroid. On the other hand, if $\varphi=0$, then $((A,
\phi),(A^*,W),0)\equiv ((A, \phi),(A^*,W))$ is a Jacobi
bialgebroid over $M$.

\begin{remark}
{\rm In \cite{jf-tj}, we proved that the double of a quasi-Jacobi
bialgebroid is a Courant-Jacobi algebroid (\cite{gm2}, \cite{jj}).
By a similar computation, we may conclude that the double of a
Jacobi-quasi bialgebroid is also a Courant-Jacobi algebroid.}
\end{remark}

The rest of this section is devoted to some important examples of
quasi-Jacobi and Jacobi-quasi bialgebroids.

\subsection{Quasi-Jacobi and Jacobi-quasi
bialgebras} A \emph{quasi-Jacobi bialgebra} is a quasi-Jacobi
bialgebroid over a point, that is a triple
$((\mathcal{G},\phi),(\mathcal{G}^*,W),\varphi)$, where
$(\mathcal{G}^*,[\cdot,\cdot]_*,W)$ is a real Lie algebra of
finite dimension with $1$-cocycle $W\in \mathcal{G}$ in its
Chevalley-Eilenberg cohomology,  $(\mathcal{G},\phi)$ is the dual
space of $\mathcal{G}^*$ endowed with a bilinear skew-symmetric
bracket $[\cdot,\cdot]$ and an element $\phi\in \mathcal{G}^*$ and
$\varphi \in \bigwedge^3\mathcal{G}^*$, such that conditions 3)-6)
of Definition \ref{def1.1} are satisfied.

By dualizing the above notion, we get a \emph{Jacobi-quasi
bialgebra}, i.e. a Jacobi-quasi bialgebroid over a point.

In the particular case where $\varphi=0$, we recover the concept
of Jacobi bialgebra \cite{im1}. When $\phi=0$ and $W=0$, we
recover the notion of quasi-Lie bialgebra due to Drinfeld
\cite{dr}.

We postpone the study of quasi-Jacobi bialgebras to a future
paper. We believe that they can be considered as the infinitesimal
invariants of Lie groups endowed with a certain type of twisted
Jacobi structures that  can be constructed from the solutions of a
twisted Yang-Baxter equation.

\subsection{The quasi-Jacobi and the Jacobi-quasi bialgebroids of a
twisted Jacobi manifold} We recall that a \emph{twisted Jacobi
manifold} \cite{jf-tj} is a differentiable manifold $M$ equipped
with a section $(\Lambda,E)$ of $\bigwedge^2(TM \times \R)$ and a
$2$-form $\omega$ such that
\begin{equation} \label{1.1}
\frac{1}{2}[(\Lambda, E), (\Lambda,
E)]^{(0,1)}=(\Lambda,E)^{\#}(\delta \omega,\omega),
\footnote{Since, for any $(\varphi,\omega) \in
\Gamma(\bigwedge^3(T^*M \times \R))$, ${\rm d}^{(0,1)}(\varphi,
\omega)=(\delta \varphi, \varphi- \delta \omega)$ and ${\rm
d}^{(0,1)}(\varphi, \omega)=(0,0)\Leftrightarrow \varphi = \delta
\omega$, equation (\ref{1.1}) means that $\frac{1}{2}[(\Lambda,
E), (\Lambda, E)]^{(0,1)}$ is the image by $(\Lambda,E)^{\#}$ of a
${\rm d}^{(0,1)}$-closed $3$-form of $TM\times \R$.}
\end{equation}
where $[\cdot,\cdot]^{(0,1)}$ denotes the Schouten bracket of the
Lie algebroid $(TM\times \R, [\cdot,\cdot], \pi)$ over $M$
modified by the 1-cocycle $(0,1)$ and $(\Lambda,E)^{\#}$ is the
natural extension of the homomorphism of
$C^{\infty}(M,\R)$-modules $(\Lambda,E)^{\#} : \Gamma(T^*M\times
\R)\to \Gamma(TM\times \R)$,
$(\Lambda,E)^{\#}(\alpha,f)=(\Lambda^{\#}(\alpha)+fE,-\langle
\alpha,E\rangle)$,  to a homomorphism from
$\Gamma(\bigwedge^k(T^*M\times \R))$ to
$\Gamma(\bigwedge^k(TM\times \R))$, $k\in \N$, given,  for any
$(\eta,\xi) \in \Gamma (\bigwedge^k( T^*M \times \R))$ and
$(\alpha_1,f_1), \ldots,(\alpha_k,f_k)\in \Gamma(T^*M \times \R)$,
by
\begin{eqnarray*}
\lefteqn{(\Lambda,E)^\#(\eta,\xi)((\alpha_1,f_1), \ldots,
(\alpha_k,f_k))}  \\
& & =(-1)^k(\eta,\xi)((\Lambda,E)^\# (\alpha_1,f_1), \cdots ,
(\Lambda,E)^\# (\alpha_k,f_k))
\end{eqnarray*}
and, for all $f\in \C$, by $(\Lambda,E)^\#(f)=f$.

 Examples of twisted
Jacobi manifolds are \emph{twisted exact Poisson manifolds} and
\emph{twisted locally conformal symplectic manifolds}, both of
them presented in \cite{jf-tj}, and also \emph{twisted contact
Jacobi manifold} described in \cite{jf-ten}.

If $(M,(\Lambda,E), \omega)$ is a twisted Jacobi manifold and $f$
a function on $M$ that never vanishes, we can define a new twisted
Jacobi structure $((\Lambda^f,E^f), \omega^f)$ on $M$, which is
said to be $f$-{\em conformal} to $((\Lambda,E), \omega)$, by
setting
$$
\Lambda^f= f \Lambda \,; \quad E^f= \Lambda^\#(\delta f)+fE \,;
\quad \omega^f = \frac{1}{f} \,\omega.
$$

In the sequel, let $(M,(\Lambda,E),\omega)$ be a twisted Jacobi
manifold and $(T^*M\times
\R,[\cdot,\cdot]_{(\Lambda,E)}^{\omega},\pi \circ
(\Lambda,E)^{\#},(-E,0))$ its canonically associated Lie algebroid
with $1$-cocycle, \cite{jf-tj}. The Lie bracket
$[\cdot,\cdot]_{(\Lambda,E)}^{\omega}$ on $\Gamma(T^*M\times \R)$
is given, for all $(\alpha,f),(\beta,g) \in \Gamma(T^*M\times
\R)$, by
\begin{eqnarray*}
[(\alpha,f),(\beta,g) ]_{(\Lambda,E)}^{\omega} & = &
[(\alpha,f),(\beta,g)]_{(\Lambda,E)} \\ &  &+\,(\delta
\omega,\omega)((\Lambda,E)^\# (\alpha,f), (\Lambda,E)^\#
(\beta,g), \cdot),
\end{eqnarray*}
where $[\cdot,\cdot]_{(\Lambda,E)}$ is the usual bracket on
$\Gamma(T^*M\times \R)$ associated to a section $(\Lambda,E)$ of
$\bigwedge^2(TM\times \R)$ (\cite{krb}, \cite{im1}):
\begin{eqnarray}\label{br-krb}
[(\alpha,f),(\beta,g)]_{(\Lambda,E)}& =& {\mathcal L
}^{(0,1)}_{(\Lambda,E)^\#(\alpha,f)}(\beta,g)- {\mathcal L
}^{(0,1)}_{(\Lambda,E)^\#(\beta,g)}(\alpha,f) \nonumber \\ & &-
\,{\rm d}^{(0,1)}((\Lambda,E)((\alpha,f),(\beta,g))).
\end{eqnarray}
We consider, on the vector bundle $TM \times \R \to \R$, the Lie
algebroid structure over $M$ with $1$-cocycle
$([\cdot,\cdot],\pi,(0,1))$ and also a new bracket
$[\cdot,\cdot]'$ on the space of its smooth sections given, for
all $(X,f),(Y,g) \in \Gamma(TM\times \R)$, by
$$
[(X,f), (Y,g)]'= [(X,f), (Y,g)]-(\Lambda,E)^\#((\delta \omega,
\omega)((X,f), (Y,g), \cdot)).
$$
We have shown in \cite{jf-tj} that the triple $((TM\times \R,
[\cdot,\cdot]',\pi,(0,1)), (T^*M\times
\R,[\cdot,\cdot]_{(\Lambda,E)}^{\omega},\pi \circ
(\Lambda,E)^{\#},(-E,0)),(\delta\omega, \omega))$ is a
quasi-Jacobi bialgebroid over $M$. Furthermore, we have:

\begin{theorem}\label{th-jq}
The triple $((TM\times \R, [\cdot,\cdot],\pi,(0,1)), (T^*M\times
\R,[\cdot,\cdot]_{(\Lambda,E)},\pi \circ
(\Lambda,E)^{\#},(-E,0)),(\Lambda,E)^{\#}(\delta\omega, \omega))$
is a Jacobi-quasi bialgebroid over $M$.
\end{theorem}
\begin{proof}
It suffices to check that the dual versions of all conditions of
Definition \ref{def1.1} are satisfied. Condition 1) can be checked
directly, using (\ref{br-krb}). For 2), we take into account that
$((\Lambda,E), \omega)$ is a twisted Jacobi structure, hence
(\ref{1.1}) holds, and we apply the general formula
\begin{eqnarray}\label{0}
(\Lambda,E)^\#([(\alpha,f),(\beta,g)]_{(\Lambda,E)})&=&[
(\Lambda,E)^\#(\alpha,f),(\Lambda,E)^\#(\beta,g)] \\
&-&
\frac{1}{2}[(\Lambda,E),(\Lambda,E)]^{(0,1)}((\alpha,f),(\beta,g),\cdot).\nonumber
\end{eqnarray}
By projection, we obtain 2). Condition 3) can be checked directly,
after a long computation. In order to prove 4), we remark that the
quasi-differential operator $d_*$ determined by $([\cdot,\cdot
]_{(\Lambda,E)},\pi \circ (\Lambda,E)^\#)$ is given \cite{im1},
for all $(R,S) \in \Gamma(\bigwedge^k(TM \times \R))$, by
$$
d_*(R,S)= ([\Lambda, R]+ k E \wedge R + \Lambda \wedge S,
-[\Lambda, R]+ (1-k) E \wedge S + [E,R]).
$$
So, $d_*(-E,0)= ([E,\Lambda],0)$, and since
$(M,(\Lambda,E),\omega)$ is a twisted Jacobi manifold,
\begin{eqnarray*}
d_*(-E,0)&=&([E,\Lambda],0)=
\frac{1}{2}[(\Lambda,E),(\Lambda,E)]^{(0,1)}((0,1),\cdot,\cdot)\\
&=& ((\Lambda,E)^\# (\delta \omega, \omega))((0,1), \cdot, \cdot).
\end{eqnarray*}
On the other hand, since $d_*^{(-E,0)}(R,S)= [(\Lambda,E),
(R,S)]^{(0,1)}$, we have
\begin{eqnarray*}
d_*^{(-E,0)}((\Lambda,E)^\# (\delta \omega, \omega))&=&
[(\Lambda,E),(\Lambda,E)^\# (\delta \omega, \omega)]^{(0,1)}  \\
&=& \frac{1}{2}
[(\Lambda,E),[(\Lambda,E),(\Lambda,E)]^{(0,1)}]^{(0,1)} =0,
\end{eqnarray*}
whence we get condition 5). Finally, 6) can be established, as in
the proof of Theorem 8.2 in \cite{jf-tj}, by a straightforward but
long computation.
\end{proof}

\vspace{2mm}

In the case of twisted Poisson manifolds, the previous results
were treated in \cite{ro1} and \cite{ks3}.

\subsection{The Jacobi-quasi bialgebroid of a quasi Jacobi manifold}\label{sec-qja}
Let $(\mathcal{G}, [\cdot,\cdot])$ be a Lie algebra, $\phi$ a
$1$-cocycle in its Chevalley-Eilenberg cohomology and
$(\cdot,\cdot)$ a nondegenerate symmetric bilinear invariant form
on $\mathcal{G}$. We denote by $\psi$  the canonical $3$-form on
$\mathcal{G}$ defined by $\psi(X,Y,Z)=\frac{1}{2}(X,[Y,Z])$, for
all $X,Y,Z \in \mathcal{G}$, and by $Q_{\psi}\in \bigwedge^3
\mathcal{G}$  its dual trivector that is given, for all $\mu, \nu,
\xi \in \mathcal{G}^*$, by
$$
Q_{\psi}(\mu, \nu, \xi)=\psi(X_{\mu},X_{\nu},X_{\xi}),
$$
where $X_{\mu},X_{\nu},X_{\xi}$ are, respectively, dual to $\mu,
\nu, \xi $ via $(\cdot,\cdot)$.

A $(\mathcal{G},\phi)$-manifold $M$ is a differentiable manifold
on which $(\mathcal{G},\phi)$ acts infinitesimally by $a^{\phi} :
\mathcal{G} \to TM\times \R$, $a^{\phi}(X)=a(X)+\langle \phi,
X\rangle$, for all $X\in \mathcal{G}$ (see Appendix). We keep the
same notation $a^{\phi}$ for the induced maps on exterior
algebras.

Let $M$ be a $(\mathcal{G},\phi)$-manifold. A section
$(\Lambda,E)\in \Gamma(\bigwedge^2 (TM\times \R))$ is said to be
{\em invariant} (under the infinitesimal action $a^\phi$) if, for
any $X\in \mathcal{G}$,
$$
{\mathcal L}^{(0,1)}_{a^\phi(X)} (\Lambda, E)=(0,0).
$$

A natural generalization of the notion of quasi Poisson manifold,
 given in \cite{aks}, is the concept of
\emph{$(\mathcal{G},\phi)$-quasi Jacobi manifold}, that we
introduce as follows.
\begin{definition}
A \emph{$(\mathcal{G},\phi)$-quasi Jacobi manifold} is a
$(\mathcal{G},\phi)$-manifold $M$ equip\-ped with an invariant
section $(\Lambda,E)\in \Gamma(\bigwedge^2 (TM\times \R))$ such
that
$$
\frac{1}{2}[(\Lambda,E),(\Lambda,E)]^{(0,1)}=a^{\phi}(Q_{\psi}).
$$
\end{definition}

A long, but not difficult computation, leads us to the following:

\begin{theorem}\label{th-jq-QJ} Let
$(M,\Lambda,E)$ be a $(\mathcal{G},\phi)$-quasi Jacobi manifold.
Then, $((TM\times \R, [\cdot,\cdot],\pi,(0,1)), (T^*M\times
\R,[\cdot,\cdot]_{(\Lambda,E)},\pi \circ (\Lambda,E)^{\#},(-E,0)),
a^{\phi}(Q_{\psi}))$ is a Jacobi-quasi bialgebroid over $M$.
\end{theorem}

\begin{remark}
\emph{If $M$ is a $\mathcal{G}$-manifold equipped with a quasi
Poisson structure, i.e an invariant bivector field $\Lambda$ on
$M$ such that $[\Lambda,\Lambda]=2a(Q_{\psi})$, a similar result
holds: \emph{The triple
$((TM,[\cdot,\cdot],id),(T^*M,[\cdot,\cdot]_{\Lambda},
\Lambda^{\#}), a(Q_{\psi}))$ is a Lie-quasi bialgebroid over $M$,
$[\cdot,\cdot]_{\Lambda}$ being the Koszul bracket associated to
$\Lambda$.}}
\end{remark}

\subsection{Triangular quasi-Jacobi and Jacobi-quasi bialgebroids}
Let $(A,[\cdot,\cdot],a,\phi)$ be a Lie algebroid with $1$-cocycle
over a differentiable manifold $M$, $\Pi$ a section of
$\bigwedge^2A$ and $Q$ a trivector on $A$ such that
\begin{equation}\label{Pi}
\frac{1}{2}[\Pi,\Pi]^\phi=Q.
\end{equation}
We shall discuss what happens on the dual vector bundle $A^*$ of
$A$ when we consider the vector bundle map $a_* : A^*\to TM$,
$a_*=a\circ \Pi^{\#}$, $\Pi^{\#} : A^*\to A$ being the bundle map
associated to $\Pi$, and the Koszul bracket $[\cdot, \cdot]_{\Pi}$
on the space $\Gamma(A^*)$ of its smooth sections given, for all
$\alpha, \beta \in \Gamma(A^*)$, by
\begin{equation} \label{koz}
[\alpha, \beta]_{\Pi}= {\mathcal L}^\phi_{\Pi^\# (\alpha)} \beta -
{\mathcal L}^\phi_{\Pi^\# (\beta)} \alpha -
d^\phi(\Pi(\alpha,\beta)).
\end{equation}
Let us set $W=-\Pi^\#(\phi)$. Taking into account that, for all
$\alpha, \beta, \gamma \in \Gamma(A^*)$,
$$
[[\alpha, \beta]_{\Pi},\gamma]_{\Pi}+c.p. =
-d^{\phi}(Q(\alpha,\beta,\gamma))-((i_{Q(\alpha,\beta,\cdot)}d^{\phi}\gamma)+c.p),
$$
we can directly prove that
\begin{theorem}
The triple $((A,[\cdot,\cdot],a,\phi),(A^*,[\cdot, \cdot]_{\Pi},
a_*, W),Q)$ is a Jacobi-quasi bialgebroid over $M$, which is
called a \emph{triangular Jacobi-quasi bialgebroid}.
\end{theorem}

Clearly, the Lie-quasi bialgebroid associated to a twisted Poisson
manifold \cite{ro1} and the Jacobi-quasi bialgebroid associated to
a twisted Jacobi manifold (see Theorem \ref{th-jq}) are special
cases of triangular Jacobi-quasi bialgebroids. Another important
type of triangular quasi-Jacobi bialgebroid is the triangular
quasi-Jacobi bialgebra, where $\Pi$ is a solution of a
Yang-Baxter's type equation.

\vspace{2mm}

Now, we consider the particular case where $Q =
\Pi^{\#}(\varphi)$, with $\varphi$ a $d^\phi$-closed $3$-form  on
$A$, and the spaces $\Gamma(A^*)$ and $\Gamma(A)$ are equipped,
respectively, with the brackets
$$
[\alpha, \beta]_{\Pi}^{\varphi}= [\alpha, \beta]_{\Pi} +
\varphi(\Pi^\# (\alpha),\Pi^\# (\beta), \cdot),   \hspace{5mm}
\mathrm{for}\;\mathrm{all}\;\; \alpha, \beta \in \Gamma(A^*),
$$
$[\cdot, \cdot]_{\Pi}$ being the Koszul bracket (\ref{koz}), and
$$
[X,Y]'=[X,Y]- \Pi^\#(\varphi(X,Y,\cdot)), \hspace{5mm}
\mathrm{for}\;\mathrm{all}\;\; X,Y\in \Gamma(A).
$$

Under the above assumptions, by a straightforward calculation, we
get:
\begin{proposition}\label{p2.12}
The vector bundle $A^*\to M$ endowed with the structure $([\cdot,
\cdot]_{\Pi}^{\varphi}, a_*)$ is a Lie algebroid over $M$ with
$1$-cocycle $W=-\Pi^\#(\phi)$.
\end{proposition}

Also, we have:

\begin{theorem}\label{th-trn}
The triple $((A,[\cdot,\cdot]',a,\phi),(A^*,[\cdot,
\cdot]_{\Pi}^\varphi, a_*, W),\varphi)$ is a triangular
quasi-Jacobi bialgebroid over $M$.
\end{theorem}
\begin{proof}
The proof is analogous to that of Theorem 8.2 in \cite{jf-tj} and
so it is omitted.
\end{proof}

\begin{remark}
\emph{Obviously, if  $A$ is $TM\times \R$ equipped with the usual
Lie algebroid structure with $1$-cocycle, $([\cdot, \cdot], \pi,
(0,1))$,  and $\Pi=(\Lambda,E)\in \Gamma(\bigwedge^2(TM \times
\R))$ satisfies (\ref{Pi}), then the manifold $M$ is endowed with
a twisted Jacobi structure. The Lie algebroid structure on
$A^*=T^*M \times \R$ given by Proposition \ref{p2.12}, is the Lie
algebroid structure canonically associated with the twisted Jacobi
structure on $M$.}
\end{remark}

\section{The structure induced on the base manifold of a
quasi-Jacobi bialgebroid}

In this section we will investigate the structure induced on the
base manifold of a quasi-Jacobi bialgebroid. Similar results hold
for a Jacobi-quasi bialgebroid.

Let $((A,\phi),(A^*,W),\varphi)$ be a quasi-Jacobi bialgebroid
over $M$. In \cite{jf-tj}, we have already considered the bracket
$\{\cdot,\cdot\}$ on $\C$ defined, for all $f,g \in \C$, by
\begin{equation} \label{4.2}
\{f,g \}=\langle d^\phi f, d_*^W g\rangle.
\end{equation}
We have proved that it is $\R$-bilinear, skew-symmetric and a
first order differential operator on each argument \cite{jf-tj}.
On the other hand, the quasi-differential operator $d$ on
$\Gamma(\bigwedge A^*)$ determined by $(a,[\cdot,\cdot])$ is a
derivation with respect to the usual product of functions.
Therefore, the map $(f,g) \mapsto \langle df, d_*g \rangle$ is a
derivation on each argument and so, there exists a bivector field
$\Lambda$ on $M$ such that, for all $f,g \in \C$,
$$
\Lambda(\delta f,\delta g)=\langle df, d_*g \rangle =- \langle dg,
d_*f \rangle.
$$
If $E$ is the vector field $a_*(\phi)=-a(W)$  on $M$ then, from
(\ref{4.2}) and because $\langle \phi, W\rangle=0$ holds
\cite{jf-tj}, we get
\begin{equation}\label{L-br}
\{f, g \}=\langle {\rm d}^{(0,1)} g, (\Lambda,E)^\#({\rm
d}^{(0,1)} f) \rangle .
\end{equation}
Since, for all $f \in \C$, $d^\phi f =(a^\phi)^*({\rm d}^{(0,1)}
f)$ and $d_*^W f= (a_*^W)^*({\rm d}^{(0,1)} f)$ \cite{jf-tj},
where $(a^\phi)^*$ and $(a_*^W)^*$ denote, respectively, the
transpose of $a^\phi$ and $a_*^W$, we obtain
\begin{equation} \label{5.4}
(\Lambda,E)^\#=- a^\phi \circ (a_*^W)^*= a_*^W \circ (a^\phi)^*.
\end{equation}
It is well known that any bracket of type (\ref{L-br}) satisfies
the following relation:
\begin{equation}\label{jac-id}
\{f,\{g,h\}\} + c.p. =
\frac{1}{2}[(\Lambda,E),(\Lambda,E)]^{(0,1)}({\rm d}^{(0,1)}
f,{\rm d}^{(0,1)} g,{\rm d}^{(0,1)} h).
\end{equation}
Therefore, for the bracket defined by (\ref{4.2}), in general, the
Jacobi identity does not hold.

\begin{proposition} \label{p5.2}
Let $((A,\phi),(A^*,W),\varphi)$ be a quasi-Jacobi bialgebroid
over $M$. Then, the bracket (\ref{4.2}) satisfies, for all $f,g,h
\in \C$, the following identity:
\begin{equation}\label{id-prop}
\{f, \{g,h\}\}+ c.p. =  a_*^W(\varphi)({\rm d}^{(0,1)} f,{\rm
d}^{(0,1)} g,{\rm d}^{(0,1)} h).
\end{equation}
In (\ref{id-prop}), $a_*^W$ denotes the natural extension of
$a_*^W: \Gamma(A^*) \to \Gamma(TM \times \R)$ to a bundle map from
$\Gamma(\bigwedge^3 A^*)$ to $\Gamma(\bigwedge^3(TM \times \R))$.
\end{proposition}
\begin{proof} Let $f,g$ and $h$ be any three functions on $\C$. Taking into account the formul{\ae} $d_*^W \{f,g \}=[d_*^Wg,
d_*^Wf]$ (see \cite{jf-tj}), (\ref{5.4}) and (\ref{0}), and the
properties of a quasi-Jacobi bialgebroid, after a simple
computation, we get
\begin{eqnarray*}
\{h, \{f, g \} \} & = & \{ \{g, f \},h \}
-\frac{1}{2}[(\Lambda,E), (\Lambda,E)]^{(0,1)}({\rm d}^{(0,1)}
f,{\rm d}^{(0,1)} g,{\rm d}^{(0,1)} h ) \\
&  & + \, \, a_*^W(\varphi)({\rm d}^{(0,1)} f,{\rm d}^{(0,1)}
g,{\rm d}^{(0,1)} h).
\end{eqnarray*}
Consequently,
\begin{equation} \label{4.4.}
\frac{1}{2}[(\Lambda,E), (\Lambda,E)]^{(0,1)}({\rm d}^{(0,1)}
f,{\rm d}^{(0,1)} g,{\rm d}^{(0,1)} h )= a_*^W(\varphi)({\rm
d}^{(0,1)} f,{\rm d}^{(0,1)} g,{\rm d}^{(0,1)} h).
\end{equation}
Hence, from (\ref{jac-id}) and (\ref{4.4.}), we obtain
(\ref{id-prop}).
\end{proof}

\vspace{2mm}

Looking at equation (\ref{4.4.}), we remark that the obstruction
for $(M,\Lambda,E)$ to be a Jacobi manifold, i.e. to have
$[(\Lambda,E), (\Lambda,E)]^{(0,1)}=(0,0)$, is the image by
$a_*^W$ of the element $\varphi$ in $\Gamma(\bigwedge^3A^*)$. This
obstruction can also be viewed as the image of $\varphi$ under the
infinitesimal action of the Lie algebroid with $1$-cocycle
$(A^*,W)$ on $M$ (see Appendix). Thus, inspired by the analogous
terms of \emph{quasi Poisson $\mathcal{G}$-manifold} (\cite{aks},
\cite{ks3}) and of \emph{$(\mathcal{G},\phi)$-quasi Jacobi
manifold} (see  Section \ref{sec-qja}), we say that the pair
$(\Lambda,E)$ defines on $M$ a $(A^*,W)$-{\em quasi Jacobi
structure}.

\vspace{2mm}

Thus, we have proved:

\begin{theorem}\label{th5.1}
Let $((A,\phi),(A^*,W),\varphi)$ be a quasi-Jacobi bialgebroid
over $M$. Then, the bracket $\{\cdot, \cdot \}: \C \times \C \to
\C$ given by
$$
\{f,g \}=\langle d^\phi f, d_*^W g\rangle,\quad for \,\,f,g \in
\C,
$$
defines a $(A^*,W)$-quasi Jacobi structure on $M$.
\end{theorem}

\begin{remark}
\emph{In the case where $((A,\phi),(A^*,W),Q)$ is a Jacobi-quasi
bialgebroid over $M$, we can easily prove that the Jacobi identity
of the bracket defined by (\ref{4.2}) is violated by the image of
$Q$ under $a^\phi$. For this reason, we shall call the structure
$(\Lambda,E)$ induced on $M$, an} $(A,\phi)$-quasi Jacobi
structure. \emph{We note that, for the proof of this result, we
use the relation $[d^{\phi}f,d^{\phi}g]_* = d^{\phi}\{f,g\}$,
$f,g\in \C$, which leads to}
\begin{equation}\label{br-L}
(\Lambda,E)^\#= a^\phi \circ (a_*^W)^*= -a_*^W \circ (a^\phi)^*.
\end{equation}
\end{remark}

\begin{examples}
\end{examples}
\vspace{-2mm}

\noindent  \emph{1) $A^*$-Quasi Poisson structures:} If
$((A,\phi),(A^*,W),\varphi)$ is a quasi-Lie bialgebroid over $M$,
i.e. both $1$-cocycles $\phi$ and $W$ are zero, Theorem
\ref{th5.1} establishes the existence of a structure on $M$,
defined by the bracket
$$
\{f,g \}=\langle df, d_* g \rangle, \hspace{8mm} f,g\in \C,
$$
on $\C$, which is associated to a bivector filed $\Lambda$ on $M$
satisfying $[\Lambda, \Lambda]=2 a_*(\varphi)$. In our
terminology, $\Lambda$ endows $M$ with a $A^*${\em -quasi Poisson
structure}. We remark that this result was obtained in \cite{ilx}
by different techniques.

\vspace{2mm}

\noindent  \emph{2) Jacobi structures:} When $\varphi=0$, i.e.
$((A,\phi),(A^*,W))$ is a Jacobi bialgebroid over $M$, the
structure $(\Lambda,E)$ on $M$ determined by Theorem \ref{th5.1}
is a Jacobi structure, and we recover the well known result of
\cite{im1}.

\vspace{2mm}

\noindent  \emph{3) Twisted Jacobi structures:} When $\varphi$ is
the image of an element $(\varphi_M, \omega_M) \in
\Gamma(\bigwedge^3 (T^*M \times \R))$ by the transpose map
$(a^\phi)^*: \Gamma(\bigwedge^3 (T^*M \times \R)) \to
\Gamma(\bigwedge^3A^*)$ of $a^\phi$, i.e.
$\varphi=(a^\phi)^*(\varphi_M, \omega_M)$, then,
$$
\frac{1}{2}[(\Lambda,E), (\Lambda,E)]^{(0,1)}  = a_*^W(\varphi) =
a_*^W((a^\phi)^*(\varphi_M, \omega_M)) \stackrel{(\ref{5.4})}{=}
(\Lambda,E)^{\#}(\varphi_M, \omega_M).
$$
Also, we have
$$
d^\phi (\underbrace{(a^\phi)^*(\varphi_M, \omega_M)}_{=
\varphi})=0 \Leftrightarrow (a^\phi)^*({\rm d}^{(0,1)}(\varphi_M,
\omega_M))=0,
$$
which means that $(\varphi_M, \omega_M)$ is ${\rm
d}^{(0,1)}$-closed on the distribution $Im(a^\phi)$. This
distribution is not, in general, involutive due to condition 2) of
Definition \ref{def1.1}. However, when $Im(a^\phi)$ is involutive,
as in the case where $a^\phi$ is surjective,
$((\Lambda,E),\omega_M)$ defines a twisted Jacobi structure on the
leaves of $Im(a^\phi)$.

\vspace{2mm}

\noindent \emph{4) The case of the quasi-Jacobi bialgebroid
associated to a twisted Jacobi manifold:} Let $(M,(\Lambda_1,E_1),
\omega)$ be a twisted Jacobi manifold and let $((TM\times \R,
[\cdot,\cdot]',\pi,$ $(0,1)), (T^*M\times
\R,[\cdot,\cdot]_{(\Lambda_1,E_1)}^{\omega},
\pi\circ(\Lambda_1,E_1)^{\#},(-E_1,0)),(\delta\omega,\omega))$ be
its associated quasi-Jacobi bialgebroid. Then, the $(T^*M\times
\R,(-E_1,0))$-quasi Jacobi structure induced on $M$ coincides with
the initial one $(\Lambda_1,E_1)$. In fact, for any $f,g \in \C$
and taking into account that ${\rm d}'^{(0,1)}f= {\rm d}^{(0,1)}
f$, ${\rm d}'$ being the quasi-differential of $TM \times \R$
determined by the structure $([\cdot,\cdot]', \pi)$, and that
$(d_*^{\omega})^{(-E_1,0)}g= - (\Lambda_1,E_1)^{\#}({\rm
d}^{(0,1)} g)$, we have
$$
\{f, g \}  =  \langle {\rm d}'^{(0,1)} f, (d_*^\omega)^{(-E_1,0)}
g \rangle  =  \langle {\rm d}^{(0,1)} f,  -
(\Lambda_1,E_1)^\#({\rm d}^{(0,1)} g) \rangle = \{f, g \}_1,
$$
where $\{\cdot,  \cdot\}_1$ denotes the bracket associated to
$(\Lambda_1,E_1)$.

Moreover, if we consider the Jacobi-quasi bialgebroid $((TM\times
\R, [\cdot,\cdot],\pi,$ $(0,1)), (T^*M\times
\R,[\cdot,\cdot]_{(\Lambda_1,E_1)},\pi\circ(\Lambda_1,E_1)^{\#},(-E_1,0)),
(\Lambda_1,E_1)^{\#}(\delta\omega,\omega))$ associated to the
twisted Jacobi manifold $(M,(\Lambda_1,E_1), \omega)$, we get that
the $(TM\times \R, (0,1))$-quasi Jacobi structure $(\Lambda,E)$
induced on $M$ is the opposite of $(\Lambda_1,E_1)$. It suffices
to remark that
\begin{eqnarray} \label{ex3.4.5}
(\Lambda,E)^{\#} & \stackrel{(\ref{br-L})}{=} & \pi^{(0,1)} \circ
((\pi \circ (\Lambda_1,E_1)^{\#})^{(-E_1,0)})^* \nonumber \\ &= &
\pi^{(0,1)} \circ ((\Lambda_1,E_1)^{\#})^* \circ (\pi^{(0,1)})^*=-
(\Lambda_1,E_1)^{\#}.
\end{eqnarray}

\vspace{2mm}

\noindent \emph{5) The induced structure on a quasi Jacobi
manifold:} We consider the Jacobi-quasi bialgebroid $((TM\times
\R, [\cdot,\cdot],\pi,(0,1)), (T^*M\times
\R,[\cdot,\cdot]_{(\Lambda_1,E_1)},\pi \circ
(\Lambda_1,E_1)^{\#},(-E_1,0)), a^{\phi}(Q_{\psi}))$ associated to
a $(\mathcal{G},\phi)$-quasi Jacobi manifold $(M,\Lambda_1,E_1)$.
Then, repeating the computation (\ref{ex3.4.5}), we conclude that,
as in the previous case, the $(TM\times \R,(0,1))$-quasi Jacobi
structure $(\Lambda,E)$ induced on $M$ is the opposite of
$(\Lambda_1,E_1)$.

\vspace{2mm}

\noindent \emph{6) The case of a triangular quasi-Jacobi
bialgebroid:} If we consider a triangular quasi-Jacobi bialgebroid
over $M$ of type $((A,[\cdot,\cdot]',a,\phi),(A^*,[\cdot,
\cdot]_{\Pi}^\varphi,$ $a_*, W),\varphi)$, presented in Theorem
\ref{th-trn}, then, for all $f \in \C$,
$$
d_*^W f = - \Pi^\#(d^\phi f)= - (\Pi^\# \circ (a^\phi)^*)({\rm
d}^{(0,1)}f).
$$
So, the bracket (\ref{4.2}) in $\C$ is given by
$$
\{f, g \}  = \langle d^\phi f, d_*^W g \rangle = \langle {\rm
d}^{(0,1)} g, (a^\phi \circ \Pi^\# \circ (a^\phi)^*){\rm
d}^{(0,1)} f \rangle.
$$
On the other hand, considering the $(A^*,W)$-quasi Jacobi
structure $(\Lambda,E)$ on $M$, we also have
$$
\{f, g\}  = \langle {\rm d}^{(0,1)} g, (\Lambda ,E)^\# ({\rm
d}^{(0,1)} f )\rangle.
$$
Hence,
$$
(\Lambda ,E)^\# =a^\phi \circ \Pi^\# \circ (a^\phi)^*,
$$
which means that $(\Lambda,E)$ is the image by $a^\phi$ of $\Pi$
and that $a^\phi$ is a type of ``twisted Jacobi morphism" between
$(A,\phi,\Pi)$ and $(TM\times \R, (0,1), (\Lambda,E))$.

\section{Quasi-Lie bialgebroids associated to quasi-Jacobi bialgebroids}

Given a Lie algebroid $(A,[\cdot,\cdot], a)$ over $M$, we can
endow the vector bundle $\tilde{A}=A \times \bkR \rightarrow M
\times \bkR$ with a Lie algebroid structure over $M \times \bkR$
as follows. The sections of $\tilde{A}$ can be identified with the
$t$-dependent sections of $A$, $t$ being the canonical coordinate
on $\bkR$, i.e., for any $\tilde{X} \in \Gamma(\tilde{A})$ and
$(x,t) \in M \times \bkR$, $\tilde{X}(x,t)=\tilde{X}_{t}(x)$,
where $\tilde{X}_{t} \in \Gamma(A)$. This identification induces,
in a natural way, a Lie bracket on $\Gamma(\tilde{A})$, also
denoted by $[\cdot,\cdot]$:
$$
[\tilde{X},\tilde{Y}](x,t)=[\tilde{X}_{t}, \tilde{Y}_{t}](x),
\hspace{5mm} \tilde{X}, \tilde{Y} \in \Gamma(\tilde{A}),\,\, (x,t)
\in M \times \R,
$$
and a bundle map, also denoted by $a$, $a: \tilde{A} \rightarrow
T(M \times \bkR)\equiv TM \oplus T \bkR$ with
$a(\tilde{X})=a(\tilde{X}_t)$, in such a way that
$(\tilde{A},[\cdot,\cdot],a)$ becomes a Lie algebroid over $M
\times \R$. If $\phi$ is a $1$-cocycle of the Lie algebroid $A$,
we know from \cite{im1} that $\tilde{A}$ can be equipped with two
other Lie algebroid structures over $M \times \bkR$, $([\cdot,
\cdot]^{\, \widetilde{ }\, \phi},\widetilde{a}^{\, \phi})$ and
$([\cdot, \cdot]^{\,\widehat{ }\,\,\phi}, \widehat{a}^{\, \phi})$
given, for all $\tilde{X}, \tilde{Y} \in \Gamma(\tilde{A})$, by
\begin{equation} \label{3.3}
[\tilde{X}, \tilde{Y}]^{\, \widetilde{ }\, \phi}= [\tilde{X}_t,
\tilde{Y}_t]+ \langle \phi,\tilde{X}_t\rangle \frac{\partial
\tilde{Y}}{\partial t}- \langle \phi,\tilde{Y}_t\rangle
\frac{\partial \tilde{X}}{\partial t};
\end{equation}
\begin{equation}\label{3.4}
\widetilde{a}^{\, \phi}(\tilde{X})=  a(\tilde{X}_t) + \langle
\phi,\tilde{X}\rangle \frac{\partial}{\partial t};
\end{equation}
and
\begin{equation} \label{3.1}
[\tilde{X}, \tilde{Y}]^{\,\widehat{ }\,\,\phi}= e^{-t}\left(
[\tilde{X}_t, \tilde{Y}_t] + \langle \phi,\tilde{X}_t\rangle
(\frac{\partial \tilde{Y}}{\partial t}- \tilde{Y})- \langle
\phi,\tilde{Y}_t\rangle (\frac{\partial \tilde{X}}{\partial t}-
\tilde{X})\right);
\end{equation}
\begin{equation} \label{3.2}
\widehat{a}^{\, \phi}(\tilde{X})=e^{-t} (\widetilde{a}^{\,
\phi}(\tilde{X})).
\end{equation}

\vspace{2mm}

Let $((A, \phi),(A^*,W),\varphi)$ be a quasi-Jacobi bialgebroid
over $M$. Then, $(A^*,[\cdot,\cdot]_*,a_*,W)$ is a Lie algebroid
with $1$-cocycle and we can consider on $\tilde{A}^*$ the Lie
algebroid structure $([\cdot, \cdot]_{*}^{\,\widehat{ }\,\,W},
\widehat{a}_*^{\, W})$ defined by (\ref{3.1}) and (\ref{3.2}).
Although $A$ is not endowed with a Lie algebroid structure, we can
still consider on $\Gamma(\tilde{A})$ a bracket $[\cdot,
\cdot]^{\,\widetilde{ }\,\,\phi}$ and a bundle map
$\widetilde{a}^{\, \phi}$ given by (\ref{3.3}) and (\ref{3.4}),
respectively. We set $\tilde{\varphi}= e^t \varphi$.

\begin{theorem} \label{th4.1}
Under the above assumptions, we have:
\begin{enumerate}
\item[1)] The triple $((A, \phi),(A^*,W),\varphi)$ is a
quasi-Jacobi bialgebroid over $M$ if and only if $(\tilde{A},
\tilde{A}^{*}, \tilde{\varphi})$ is a quasi-Lie bialgebroid over
$M \times \R$. \item[2)] If $\tilde{\Lambda}$ is the induced
$\tilde{A}^*$-quasi Poisson structure on $M \times \R$, then it is
the ``quasi Poissonization" of the induced $(A^*,W)$-quasi Jacobi
structure $(\Lambda,E)$ on $M$.
\end{enumerate}
\end{theorem}
\begin{proof}
1) Let us suppose that $((A, \phi),(A^*,W),\varphi)$ is a
quasi-Jacobi bialgebroid over $M$ and let $\tilde{X}$, $\tilde{Y}$
and $\tilde{Z}$ be three arbitrary sections in $\Gamma(\tilde{A})$
and $\tilde{f} \in \CR$. A straightforward computation gives
\begin{equation} \label{3.5}
[\tilde{X}, \tilde{f} \tilde{Y}]^{\, \widetilde{ }\, \phi}=
\tilde{f} [\tilde{X}, \tilde{Y}]^{\, \widetilde{ }\, \phi} +
(\widetilde{a}^{\, \phi}(\tilde{X})\tilde{f})\tilde{Y}.
\end{equation}
\noindent Moreover,
\begin{eqnarray*}
\lefteqn{\widetilde{a}^{\, \phi}([\tilde{X},  \tilde{Y}]^{\,
\widetilde{ }\, \phi}) = [a(\tilde{X}_t), a(\tilde{Y}_t)] -
a_*(\varphi(\tilde{X}_t,\tilde{Y}_t, \cdot))+\langle
\phi,[\tilde{X}_t,\tilde{Y}_t] \rangle
\frac{\partial}{\partial t}} \nonumber \\
& & + \, \langle \phi,\tilde{X} \rangle \left(a(\frac{\partial
\tilde{Y}}{\partial t})+\langle \phi,\frac{\partial
\tilde{Y}}{\partial t} \rangle \frac{\partial }{\partial t}\right)
-\langle \phi,\tilde{Y} \rangle \left(a(\frac{\partial
\tilde{X}}{\partial t})+\langle \phi,\frac{\partial
\tilde{X}}{\partial t} \rangle \frac{\partial }{\partial t}
\right)
\nonumber \\
&=&[\widetilde{a}^{\, \phi}(\tilde{X}),\widetilde{a}^{\,
\phi}(\tilde{Y})] - \widehat{a}_{*}^{
W}(\tilde{\varphi}(\tilde{X},\tilde{Y}, \cdot)),
\end{eqnarray*}
where, in the last equality, we have used $d \phi- \varphi (W,
\cdot, \cdot)=0$. On the other hand,
\begin{eqnarray}\label{3.7}
\lefteqn{[[\tilde{X},  \tilde{Y}]^{\, \widetilde{ }\, \phi},
\tilde{Z}]^{\, \widetilde{ }\, \phi} + c.p.  =
\left([[\tilde{X}_t, \tilde{Y}_t], \tilde{Z}_t] - d
\phi(\tilde{X}_t, \tilde{Y}_t)\frac{\partial \tilde{Z}_t}{\partial
t}\right)+ c.p.} \nonumber
\\
& = & -d_{*}^{W}(\varphi(\tilde{X}_t,\tilde{Y}_t,\tilde{Z}_t))-
\left(\left(i_{\varphi(\tilde{X}_t,\tilde{Y}_t,\cdot)}
d_{*}^{W}\tilde{Z}_t
 +\varphi(W,\tilde{X}_t,\tilde{Y}_t)\frac{\partial
\tilde{Z}_t}{\partial t} \right) +c.p.\right) \nonumber \\
& = & -\widehat{d}_{*}^{\,\,
W}(\tilde{\varphi}(\tilde{X},\tilde{Y},\tilde{Z}))-
\left(i_{\tilde{\varphi}(\tilde{X},\tilde{Y},\cdot)}\widehat{d}_{*}^{\,\,
W}\tilde{Z}+ c.p.\right),
\end{eqnarray}
where $\widehat{d}_{*}^{\,\, W}$ is the differential operator of
$(\tilde{A}^*,[\cdot, \cdot]_{*}^{\,\widehat{
}\,\,W},\widehat{a}_{*}^{ W})$. Because $d^{\phi} \varphi=0$, we
get
\begin{equation} \label{3.8}
\widetilde{d}^{\, \phi}\tilde{\varphi}=0,
\end{equation}
where $\widetilde{d}^{\, \phi}$ is the quasi-differential operator
determined by the structure $([\cdot,\cdot]^{\, \widetilde{ }\,
\phi}, \widetilde{a}^{\, \phi} )$ on $\tilde{A}$. Finally, after a
very long computation we obtain
\begin{equation} \label{3.9}
\widehat{d}_{*}^{\, \,W}[\tilde{P},\tilde{Q}]^{\, \widetilde{ }\,
\phi}= [\widehat{d}_{*}^{\, \,W}\tilde{P},\tilde{Q}]^{\,
\widetilde{ }\, \phi} + (-1)^{p+1}[\tilde{P},\widehat{d}_{*}^{\,\,
W}\tilde{Q}]^{\, \widetilde{ }\, \phi},
\end{equation}
for $\tilde{P} \in \Gamma(\bigwedge^p \tilde{A})$ and $\tilde{Q}
\in \Gamma(\bigwedge \tilde{A})$. From relations (\ref{3.5}) to
(\ref{3.9}), we conclude that $(\tilde{A}, \tilde{A}^{*},
\tilde{\varphi})$ is a quasi-Lie bialgebroid over $M \times \R$.

Now, let us suppose that  $(\tilde{A}, \tilde{A}^{*},
\tilde{\varphi})$ is a quasi-Lie bialgebroid over $M \times \R$
and take three sections $X$, $Y$ and $Z$ of $A$ and $f \in \C$.
These sections can be viewed as sections of $\tilde{A}$ that don't
depend on $t$, as well as the function $f$ can also be viewed as a
function on $\CR$. Condition 1) of Definition \ref{def1.1} is
immediate from $[X, fY]^{\, \widetilde{ }\, \phi}= f [X, Y]^{\,
\widetilde{ }\, \phi}+ (\widetilde{a}^{\, \phi}(X)f)Y$. The
condition $\widetilde{a}^{\, \phi}([X,Y]^{\, \widetilde{ }\,
\phi})=[\widetilde{a}^{\, \phi}(X),\widetilde{a}^{\, \phi}(Y)] -
\widehat{a}_{*}^{ W}(\tilde{\varphi}(X,Y, \cdot))$ is equivalent
to conditions 2) and 4) of Definition \ref{def1.1}. From
$\widetilde{d}^{\, \phi}\tilde{\varphi}=0$ we deduce
$d^{\phi}\varphi=0$. Finally, by similar computations, we obtain
the two remaining conditions that lead to the conclusion that
$((A, \phi),(A^*,W),\varphi)$ is a quasi-Jacobi bialgebroid over
$M$.

\vspace{2mm}

\noindent 2) Let $\tilde{\Lambda}$ be the $\tilde{A}^*$-quasi
Poisson structure induced by $(\tilde{A}, \tilde{A}^{*},
\tilde{\varphi})$ on $M \times \R$. For all $\tilde{f}, \tilde{g}
\in \CR$, we have
$$
\{\tilde{f}, \tilde{g} \}= \tilde{\Lambda}(\delta \tilde{f},\delta
\tilde{g})
$$
and, on the other hand,
$$
\{\tilde{f}, \tilde{g} \}  =  \langle \tilde{d}^{\,\phi}\tilde{f},
\hat{d}_*^{\,W} \tilde{g} \rangle = e^{-t} ( \langle d \tilde{f},
d_* \tilde{g} \rangle +\frac{\partial \tilde{f}}{\partial
t}a_*(\phi)(\tilde{g})+\frac{\partial \tilde{g}}{\partial
t}a(W)(\tilde{f})).
$$
If $(\Lambda,E)$ is the $(A^*,W)$-quasi Jacobi structure induced
by $((A, \phi),(A^*,W),\varphi)$ on $M$, since $E=a_*(\phi)=-a(W)$
and $\Lambda(\delta \tilde{f},\delta \tilde{g})=\langle d
\tilde{f}, d_* \tilde{g} \rangle$, we get that
$\tilde{\Lambda}=e^{-t}(\Lambda + \frac{\partial}{\partial
t}\wedge E)$.
\end{proof}

\vspace{3mm}

For the case of Jacobi-quasi bialgebroids we can prove a similar
result. Let $((A,\phi),(A^*,W),Q)$ be a Jacobi-quasi bialgebroid
over $M$. We consider on $\tilde{A}$ the Lie algebroid structure
$([\cdot, \cdot]^{\,\widehat{ }\,\,\phi},\widehat{a}^{\, \phi})$
defined by (\ref{3.1}) and (\ref{3.2}),  on $\tilde{A}^*$ the
structure $([\cdot, \cdot]_{*}^{\,\widetilde{ }\,\,W},
\widetilde{a}_*^{\, W})$ defined by (\ref{3.3}) and (\ref{3.4}),
and we set $\tilde{Q}=e^tQ$.

\begin{theorem}
Under the above assumptions, we have:
\begin{enumerate}
\item[1)] The triple $((A, \phi),(A^*,W),Q)$ is a Jacobi-quasi
bialgebroid over $M$ if and only if $(\tilde{A}, \tilde{A}^{*},
\tilde{Q})$ is a Lie-quasi bialgebroid over $M \times \R$.
\item[2)] If $\tilde{\Lambda}$ is the  $\tilde{A}$-quasi Poisson
structure induced on $M \times \R$, then it is the ``quasi
Poissonization" of the $(A,\phi)$-quasi Jacobi structure
$(\Lambda,E)$ induced on $M$.
\end{enumerate}
\end{theorem}

\section{Appendix: Action of a Lie algebroid with $1$-cocy\-cle}
In this Appendix, we extend the definition of Lie algebroid action
\cite{ksm} to that of Lie algebroid with $1$-cocycle action.

Let $(A,[\cdot,\cdot],a)$ be a Lie algebroid on $M$ and $\varpi :
F\to M$ a fibered manifold with base $M$, i.e. $\varpi : F\to M$
is a surjective submersion onto $M$. We recall that an
 \emph{infinitesimal action of
$A$ on $F$} (\cite{ksm}) is a $\R$-linear map
$\mathbf{ac} : \Gamma(A) \to \Gamma(TF)$ such that: \\
(1) for each $X\in \Gamma(A)$, $\mathbf{ac}(X)$ is projectable to
$a(X)$, \\
(2) the map $\mathbf{ac}$ preserves brackets, \\
(3) $ \mathbf{ac}(fX) = (f\circ \varpi)\mathbf{ac}(X) $, for all
$f\in \C$ and  $X\in \Gamma(A)$.

\begin{definition}\label{def-ac}
 An \emph{infinitesimal action of $(A,\phi)$ on $F$} is a
$\R$-linear map $\mathbf{ac^{\phi}} : \Gamma(A) \to
\Gamma(TF\times \R)$ given, for each $X\in \Gamma(A)$, by
$$
\mathbf{ac^{\phi}}(X)= \mathbf{ac}(X)+\langle \phi,X\rangle,
$$
where $\mathbf{ac}$ is an infinitesimal action of $A$ on $F$.
\end{definition}

In the particular case where $M$ is a point and therefore $A$ is a
Lie algebra, we obtain, from Definition \ref{def-ac}, the notion
of infinitesimal action of a Lie algebra with $1$-cocycle on a
manifold $F$, used in the definition of quasi Jacobi structures,
in section \ref{sec-qja}.

If, in the  Definition \ref{def-ac}, $F=M$ and $\varpi : M\to M$
is the identity, we get the concept of infinitesimal action of
$(A,\phi)$ on the base manifold $M$ that we have used to
characterize the structure induced on the base manifold of a
quasi-Jacobi bialgebroid.

\vspace{.5cm}

\noindent {\bf Acknowledgments} \, Research of J. M. Nunes da
Costa supported by CMUC-FCT and POCI/MAT/58452/2004.

\end{document}